\documentclass[11pt]{article}
\usepackage{amssymb,latexsym,amsmath}
\title{Energy decay for damped wave equations on partially rectangular domains}
\author{Nicolas Burq \footnote{ D\'epartement de Math\'ematiques, Universit\'e Paris Sud-Orsay,
91405 Orsay Cedex, France, and Institut universitaire de France.
nicolas.burq@math.u-psud.fr }\and
Michael Hitrik\footnote{Department of Mathematics,  University of
California, Los Angeles,  CA 90095-1555,
USA. hitrik@math.ucla.edu}}
\date{}

\def\wrtext#1{\relax\ifmmode{\leavevmode\hbox{#1}}\else{#1}\fi}
\def\abs#1{\left|#1\right|}
\def\begeq{\begin{equation}}
\def\endeq{\end{equation}}

\def\neigh{neighborhood}

\def\Re{{\rm Re\,}}
\def\Im{{\rm Im\,}}

\newcommand{\eps}{\epsilon}
\def\part#1{\frac{\partial}{\partial #1}}

\def\norm#1{||\,#1\,||}

\newcommand{\real}{\mbox{\bf R}}
\newcommand{\comp}{\mbox{\bf C}}

\newcommand{\nat}{\mbox{\bf N}}

\renewcommand{\Re}{\mbox{\rm Re\,}}
\renewcommand{\Im}{\mbox{\rm Im\,}}

\renewcommand{\exp}{\mbox{\rm exp\,}}

\newtheorem{dref}{Definition}[section]
\newtheorem{lemma}[dref]{Lemma}
\newtheorem{theo}[dref]{Theorem}
\newtheorem{prop}[dref]{Proposition}

\newenvironment{proof}{\vspace{.3cm}\noindent{{\em Proof:}}}{\hfill$\Box$}

\begin{document}
\maketitle

\begin{abstract}
We consider the wave equation with a damping term on a partially rec\-tan\-gular planar
domain, assuming that the damping is concentrated close to the
non-rectangular part of the domain. Polynomial decay estimates for the
ener\-gy of the solution are established.
\end{abstract}


\vskip 2mm \noindent {\bf Mathematics Subject Clas\-si\-fica\-tion
2000:} 35L05, 47A10, 47B44, 47D06, 49J20


\section{Introduction and statement of result}\label{section0}
\setcounter{equation}{0} \setcounter{dref}{0}

The purpose of this note is to show how the methods of~\cite{BZ1},~\cite{BZ2} apply to estimate the energy
decay rates for the damped wave equation on a class of planar domains, including some ergodic billiards. In a situation
when the geometric control condition of~Bardos, Lebeau, and Rauch~\cite{BLR} does not hold,
we obtain a polynomial decay estimate for the
energy of the damped wave, with respect to a stronger norm of the initial data.
In order to formulate the main result, we shall begin by recalling
some standard notation and assumptions.

Let $\Omega\subset \real^n$, $n\geq 2$, be a bounded connected domain, with
$\partial \Omega\in C^{\infty}$. When $0\leq a\in
L^{\infty}(\Omega)$ is a non-negative smooth function on $\Omega$, we
consider the following initial-boundary value problem, \begeq
\label{eq1} \left\{ \begin{array}{ll}
\left(-D_t^2-\Delta+2ia(x)D_t\right)u=0,\:\:(t,x)\in
\real_+\times \Omega, \\
u=0\quad \wrtext{on}\quad \real_{+}\times \partial \Omega, \\
u |_{t=0}= u_0\in H_0^1(\Omega),\quad D_t u|_{t=0}=u_1\in L^2(\Omega).
\end{array} \right.
\endeq
Here $D_t=\partial_t/i$ and we shall assume throughout that the damping coefficient $a$ does not vanish identically.

Associated with the evolution problem (\ref{eq1}) is the solution
operator ${\cal U}(t)=e ^{it {\cal A}}$, $t\geq 0$, acting in the
Hilbert space of the Cauchy data ${\cal H}=H_0^1(\Omega)\times L^2(\Omega)$ and mapping
$(u_0,u_1)\in {\cal H}$ to $(u(t,\cdot),D_tu(t,\cdot))$. Here we equip
${\cal H}$ with the norm
$$
\left |\left|\,\left(\begin{array}{cc}
u_0 \\
u_1
\end{array}\right)\,\right |\right|_{{\cal H}}^2=\norm{\nabla
u_0}_{L^2}^2+\norm{u_1}_{L^2}^2.
$$
The infinitesimal generator ${\cal A}$ of the semigroup ${\cal U}(t)$ is the operator
\begeq
\label{eq1.2}
{\cal A} = \left( \begin{array}{cc}
   0 & 1 \\
  -\Delta & 2ia(x)
\end{array} \right): {\cal H}\rightarrow {\cal H},
\endeq
with the domain $ D({\cal A})=\left(H^2\cap H^1_0\right)\times H_0^1$. It follows that the
spectrum of ${\cal A}$ is discrete, and
from~\cite{Lebeau} we recall that if $\lambda\in \comp$ is an
eigenvalue of ${\cal A}$ then $0<\Im \lambda \leq 2\norm{a}_{L^{\infty}}$.

For future reference, let us notice that when $U\in D({\cal A})$, then $\Im \langle{{\cal A}U,U\rangle}_{{\cal H}}\geq 0$
and hence
\begeq
\label{eq1.5}
\norm{\left(\lambda-{\cal A}\right)^{-1}}_{{\cal L}({\cal H}, {\cal H})} \leq \frac{1}{\abs{\Im \lambda}},\quad
\Im \lambda<0.
\endeq

\medskip
In this note we shall be concerned with the energy of the solution $u(x,t)$ of (\ref{eq1}) at time $t$,
\begeq
\label{eq2}
E(u,t)=\frac{1}{2}\int_{\Omega}\left(\abs{\nabla_x u}^2+\abs{D_t u}^2\right)\,dx.
\endeq
It is easily seen that $E(u,t)$ is nonincreasing as $t\rightarrow \infty$ and from~\cite{Lebeau} we may also recall that
$E(u,t)\rightarrow 0$ as $t\rightarrow \infty$, for each $(u_0,u_1)\in {\cal H}$.
Under the geometric control condition stating that there exists $T_0>0$ such any billiard trajectory in $\Omega$
of length $\geq T_0$ meets the open set $\{x; a(x)>0\}$, the uniform exponential decay of the energy has
been established by Bardos, Lebeau, and Rauch~\cite{BLR}. In the general case, without any assumptions on the underlying
dynamics, it has been proved by Lebeau~\cite{Lebeau} that the decay
rate of the energy is always logarithmic, provided that the initial
data in (\ref{eq1}) are measured with respect to a stronger norm. In
this note, we shall derive a polynomial decay
estimate for the energy, for a class of planar domains and damping
regions, in the case when the geometric control condition of~\cite{BLR} fails to hold.

We shall consider the class of partially rectangular domains
$\Omega\subset \real^2$. By this we mean that
$\Omega$ is connected, has a boundary that is piecewise $C^{\infty}$, and contains a rectangle $R\subset \Omega$,
such that if we decompose the boundary of $R$ into pairs of parallel segments, $\partial R=\Gamma_1\cup \Gamma_2$, then
$\Gamma_j\subset \partial \Omega$ for at least one $j$, say for $j=1$. We shall write $\Omega=R\cup W$, where $W$ is
the non-rectangular part of $\Omega$.

\vskip 2mm
\noindent
{\it Example}. The Bunimovich stadium $S$, defined as the union of a rectangle $R=\{(x,y); x\in [0,1],\,
y\in [0,2\beta]\}$, $\beta>0$, with two semicircular regions centered at $(0,\beta)$ and $(1,\beta)$, with
radius $\beta$, which lie outside of $R$, is a partially rectangular domain. It is well-known~\cite{Bun}
that the geodesic flow in $S$, obeying the law of reflection at the boundary, is ergodic.
In what follows, the ergodicity of the underlying classical flow in $\Omega$ will not play any special role in our
considerations.

\begin{theo}
Let $\Omega\subset \real^2$ be a partially rectangular domain, $\Omega=R\cup W$, and let
$0\leq a\in C(\overline{\Omega})$ be such that $a>0$ in $\overline{W}$, the closure of $W$ in $\Omega$. Then for each
$k>0$ there exists a constant $C_k>0$ such that for each $(u_0,u_1)\in D({\cal A}^k)$ we have
\begeq
\label{eq3}
E(u,t)^{1/2}\leq C_k \frac{(\log t)^{\frac{k}{2}+1}}{t^{\frac{k}{2}}}\norm{(u_0,u_1)}_{D({\cal A}^k)},\quad t\geq 2.
\endeq
\end{theo}

\bigskip
\noindent
{\it Remark}. We notice that in the case when the damping coefficient $a$ vanishes outside of a
\neigh{} of $\overline{W}$ in $\Omega$, the geometric control condition of~\cite{BLR} fails, due to the
existence of the invariant set for the classical flow constituted by the bouncing ball orbits
parallel to the pair of segments $\Gamma_2$ and staying within the
rectangle $R$. Therefore (see~\cite{BLR}), no
uniform decay estimate for the energy with respect to the energy norm of the initial data is possible.

\medskip
\noindent
Theorem 1.1 is proved in section 2. We shall finish this section by
briefly discussing the issue of optimality of the bound
(\ref{eq3}). First, arguing heuristically on the classical level, let us notice that every
ray starting in the interior of the rectangle and forming an angle of
size $\frac{1}{k}$, $k\gg 1$, with the direction of $\Gamma_2$ will reach the non-rectangular part of the
domain after a time which is of the order $\sim k$. When discussing
the quantum picture, for simplicity, we shall replace the rectangle $R$ by an infinite strip of the form $\real_x\times
[0,\pi]_y$ and consider the wave evolution of the functions
$e_k(x,y)=\varphi(x) \sin ky$, $\varphi\in {\cal S}(\real)$,
$\widehat{\varphi}\in C^{\infty}_0(\real)$,
$\norm{\varphi}_{L^2}=1$. Here $\widehat{\varphi}(\xi)=\int
e^{-ix\xi}\varphi(x)\,dx$ is the Fourier transform of $\varphi$. Then $e_k$ is an ${\cal O}(1)$-quasimode for the
Dirichlet Laplacian $\Delta_D$ in $\real_x\times [0,\pi]_y$, and a simple calculation shows that
\begeq
\label{eq3.5}
\left(\cos t\sqrt{-\Delta_D} \right)e_k=\left(\cos tk\right) e_k+{\cal
O}\left(\frac{t}{k}\right),\quad \abs{t}\leq k\rightarrow \infty.
\endeq
These heuristic considerations suggest that the optimal energy decay estimate in
Theorem 1.1, modulo an $\varepsilon$-loss, should be of the form
\begeq
\label{eq3.6}
E(u,t)^{1/2}\leq \frac{{\cal
O}_{\varepsilon}(1)}{t^{1-\varepsilon}}\norm{(u_0,u_1)}_{D({\cal
A})},\quad \varepsilon > 0.
\endeq
(Remark that taking $t= \varepsilon k$, $\varepsilon \ll 1$
in~\eqref{eq3.5} and letting $k$ tend to infinity shows that
\eqref{eq3.6} cannot hold if $\varepsilon <0$.)\par In section 3
we shall describe a class of examples of $C^{\infty}$-damping
coefficients $a$ for which we are able to improve the result of
Theorem 1.1 and obtain an estimate of the form (\ref{eq3.6}).

\begin{theo}
Assume that $a\in C^{\infty}(\overline{\Omega})$ is a damping coefficient
satisfying the assumptions of the beginning of section{\rm ~\ref{section3}} and in
particular {\rm (\ref{eq57.8})}, {\rm (\ref{eq58})}, for some $m>4$. Then we have
for $k>0$ and $\varepsilon>0$,
\begeq
\label{eq72}
E(u,t)^{1/2}\leq \frac{{\cal O}_{k,\varepsilon}(1)}{t^{\frac{k}{1+\frac{4}{m}}-\varepsilon}}
\norm{(u_0,u_1)}_{D({\cal A}^k)},\quad t\geq 2.
\endeq
In particular if~\eqref{eq58} is satisfied for all $m>4$, then for
each $k>0$ and each $\varepsilon>0$ there exists
$C_{k,\varepsilon}>0$ such that for each $(u_0,u_1)\in D({\cal
A}^k)$ we have
\begeq
\label{eq71}
E(u,t)^{1/2}\leq \frac{C_{k,\varepsilon}}{t^{k-\varepsilon}}\norm{(u_0,u_1)}_{D({\cal
A}^k)},\quad t\geq 2.
\endeq
\end{theo}

\medskip
\noindent
At the present time, we do
not know whether one can obtain an estimate of the form
(\ref{eq3.6}) for a general continuous, or even smooth, damping term,
concentrated near the non-rectangular part of the domain.

\medskip
\noindent
{\it Remark}. When this paper was complete, the authors learned of a very recent work by Kim Dang Phung~\cite{Phung},
where a polynomial decay estimate for the energy, with an unspecified decay rate, was obtained on a partially
cubic domain, with the damping acting in a \neigh{} of the boundary except between a pair of parallel square faces
of the cube. An initial glance at~\cite{Phung} shows that the methods of the present paper are completely different.

\medskip
\noindent {\bf Acknowledgment}. The second author is grateful to
the National Science Foundation for partial support under the
grant DMS-0304970. He would also like to thank Rowan Killip for a
stimulating discussion.

\section{Proof of Theorem 1.1}\label{section1} \setcounter{equation}{0}
Following~~\cite{BZ1} and \cite{Lebeau}, we shall use the stationary methods, and the main point will be to estimate
the resolvent $(\lambda-{\cal A})^{-1}: {\cal H}\rightarrow {\cal H}$ for $\lambda \in \real$, $\abs{\lambda}\gg 1$.

\medskip
\noindent
The following result, closely related to Proposition 6.1 of~\cite{BZ1}, is our starting point.

\begin{prop}
Let $R=[0,1]_x\times [0,\pi]_y\subset \real^2$ be a rectangle, and let us consider the stationary problem
\begeq
\label{eq4}
\left(-\Delta+2ia\lambda-\lambda^2\right)u=f+\partial_x g,\quad u|_{\partial R}=0.
\endeq
Here $\lambda \in \real$, $\abs{\lambda}\geq 1$ and $f$, $g \in L^2(R)$. Assume that $0\leq a\in C(R)$
is such that
\begeq
\label{eq4.5}
[0,1]\ni x\mapsto \inf_{y\in [0,\pi]} a(x,y)\,\wrtext{does not vanish identically}.
\endeq
Then we have
\begeq
\label{eq4.6} \norm{u}_{L^2(R)}^2\leq {\cal O}(1)\left(\norm{f}_{L^2(R)}^2+\norm{g}_{L^2(R)}^2+
\lambda^2\int_R a(x,y)\abs{u}^2\,dx\,dy\right).
\endeq
\end{prop}
\begin{proof}
When establishing (\ref{eq4.6}), it suffices to do so when $a=a(x)$
is a function of $x$ only. Indeed, assume that (\ref{eq4.6}) has
already been established in this special case. If $a\in C(R)$ satisfies (\ref{eq4.5}), we can take
$0\leq a_1\in C([0,1])$, not identically zero and such that $0\leq a_1(x)\leq a(x,y)$ when $(x,y)\in R$. Now if
$u\in H^1_0(R)$ satisfies
\begeq \label{eq5}
\left(-\Delta+2i\lambda a-\lambda^2\right)u=f+\partial_x g,
\endeq
then
\begeq
\label{eq6}
\left(-\Delta+2i\lambda a_1-\lambda^2\right)u=f+\partial_x
g+2i\lambda(a_1-a)u.
\endeq
Applying (\ref{eq4.6}) to (\ref{eq6}), we get
\begeq
\label{eq7}
\norm{u}_{L^2}^2\leq {\cal O}(1) \left(\norm{f}_{L^2}^2+\norm{g}_{L^2}^2+\lambda^2\int_R
(a-a_1)\abs{u}^2\,dx\,dy+\lambda^2 \int_R a_1 \abs{u}^2\,dx\,dy \right),
\endeq
and the bound (\ref{eq4.6}) follows in the general case. In what
follows, when proving (\ref{eq4.6}), we shall therefore assume that
$0\leq a=a(x)\in C([0,1])$ is a function of $x$ only, which does not vanish identically.

When analyzing (\ref{eq4}), we follow~\cite{BZ1},~\cite{BZ2} and separate variables.
The Dirichlet realization of $-\partial_y^2$ on $L^2((0,\pi))$
has the eigenvalues $k^2$, $k=1,2,\ldots$ with the corresponding
eigenfunctions $e_k(y)=\sqrt{2/\pi}\sin ky$, forming an orthonormal basis
in $L^2((0,\pi))$. Writing
$$
u(x,y)=\sum_{k=1}^{\infty} u_k(x) e_k(y),
$$
and similarly for $f$ and $g$,
$$
f(x,y)=\sum_{k=1}^{\infty} f_k(x) e_k(y),\quad g(x,y)=\sum_{k=1}^{\infty} g_k(x) e_k(y),
$$
we see that $u_k$, $k=1,2,\ldots$ satisfy
\begeq
\label{eq8}
\left(-\partial_x^2+2ia(x)\lambda+k^2-\lambda^2\right)u_k(x)=f_k(x)+\partial_x
g_k(x),\quad u_k(0)=u_k(1)=0.
\endeq

When analyzing (\ref{eq8}), we shall first consider the case when $k\in \nat$ is such that $k \geq \abs{\lambda}$.
We then claim that
\begeq
\label{eq9}
\norm{u_k}_{L^2}^2\leq {\cal O}(1)\left(\norm{f_k}_{L^2}^2+\norm{g_k}_{L^2}^2\right),
\endeq
and here we do not need the assumption that $a$ does not vanish identically.

Indeed, multiplying the first equation in (\ref{eq8}) by $\overline{u_k}$, we obtain, integrating by parts
and taking the real part, that
\begin{eqnarray}
\label{eq10} & &
\int_0^1 \abs{u'_k(x)}^2\,dx+(k^2-\lambda^2)\int_0^1 \abs{u_k(x)}^2\,dx
\\ \nonumber
& &= \Re\int_0^1\left(f_k(x)\overline{u_k}(x)-g_k(x)\partial_x\overline{u_k}(x)\right)\,dx.
\end{eqnarray}
Here the right hand side does not exceed
$$
\left(\norm{f_k}_{L^2}+\norm{g_k}_{L^2}\right)\left(\norm{u_k}_{L^2}+\norm{u'_k}_{L^2}\right)
\leq
2\left(\norm{f_k}_{L^2}+\norm{g_k}_{L^2}\right)\norm{u'_k}_{L^2},
$$
where we also used the Poincar\'e inequality
$$
\norm{u_k}_{L^2}\leq \norm{u_k'}_{L^2},\quad u_k\in H^1_0((0,1)).
$$
Another application of the Poincar\'e inequality combined with
(\ref{eq10}) then immediately gives (\ref{eq9}).

We now come to discuss the second case when $k< \abs{\lambda}$. Let $\tau\in \real$ be such that $\tau^2=\lambda^2-k^2$.
Then we have
\begeq
\label{eq11}
\left(-\partial_x^2+2ia(x)\tau-\tau^2\right)u_k(x)=f_k(x)+\partial_x g_k(x)+
2i(\tau-\lambda)a(x)u_k(x),
\endeq
$$
u_k(0)=u_k(1)=0.
$$
Now the operator in the left hand side of (\ref{eq11}) is a one-dimensional stationary damped wave operator, for which
it is well-known~\cite{CoxZuazua} that the resolvent bound
\begeq
\label{eq12}
R_0(\tau):=\left(-\partial_x^2+2ia\tau-\tau^2\right)^{-1}={\cal O}\left(\frac{1}{1+\abs{\tau}}\right):
L^2((0,1))\rightarrow L^2((0,1))
\endeq
holds true, since $0\leq a\in C([0,1])$ is not identically zero. See also~\cite{Hi}. It follows that
\begeq
\label{eq13}
R_0(\tau)={\cal O}(1): L^2((0,1))\rightarrow H^1_0((0,1)),\quad \tau\in \real,
\endeq
and hence by duality, the same bound holds when $R_0(\tau)$ is viewed as an operator from  $H^{-1}((0,1))$ to $L^2((0,1))$. Applying these observations to (\ref{eq11}), we get
\begeq
\label{eq14}
\norm{u_k}_{L^2}^2\leq {\cal O}(1) \left(\norm{f_k}_{L^2}^2+\norm{g_k}_{L^2}^2\right)+
{\cal O}(1)\lambda^2 \norm{a^{1/2}u_k}_{L^2}^2.
\endeq
Combining the bounds (\ref{eq9}) and (\ref{eq14}) and summing with respect to $k$, we get the
proposition.
\end{proof}

\begin{prop}
Let $\Omega\subset \real$ be a partially rectangular domain, $\Omega=R\cup W$, and let $0\leq a\in C(\overline{\Omega})$
be such that $a>0$ in $\overline{W}$. Then if we put $R(\lambda)=(-\Delta+2ia\lambda-\lambda^2)^{-1}$, $\lambda\in
\real$, then we have
\begeq
\label{eq15}
R(\lambda)={\cal O}(\abs{\lambda}): L^2(\Omega)\rightarrow L^2(\Omega),\quad \abs{\lambda}\geq 1.
\endeq
\end{prop}
\begin{proof}
We may assume that $R=[0,1]_x\times [0,\pi]_y$, with the sides
parallel to the $x$-axis contained in the boundary of $\Omega$.

\noindent
When $f\in L^2(\Omega)$, let $u\in H^1_0(\Omega)$
be the solution of
\begeq
\label{eq16}
\left(-\Delta+2ia \lambda-\lambda^2\right)u=f,\quad u|_{\partial \Omega}=0,\quad \abs{\lambda}\geq 1.
\endeq
We now let $0\leq \chi\in C^{\infty}_0((0,1))$, $0\leq \chi \leq 1$, be a cut-off function such
that $\chi=1$ on $[\eps,1-\eps]$. Here $\eps>0$ small is to be
chosen. Then it follows that $\chi u$, viewed as a function on
$R$, vanishes on its boundary and satisfies in the interior
\begeq
\label{eq16.5}
\left(-\Delta+2ia\lambda-\lambda^2\right)\chi u=\chi
f+[-\Delta,\chi]u.
\endeq
Here $[-\Delta,\chi]u=\chi''(x)u-2\partial_x (\chi'(x)u)$. An
application of Proposition 2.1 gives therefore that
\begeq
\label{eq17}
\norm{\chi u}_{L^2}^2\leq {\cal
O}(1)\left(\norm{f}_{L^2}^2+\int_{\omega_{\eps}}
\abs{u}^2\,dx\,dy+\abs{\lambda}^2 \int_R a(x,y) \abs{u}^2
dx\,dy\right).
\endeq
Here $\omega_{\eps}$ is a \neigh{} of the support of $\chi'(x)$,
and choosing $\eps>0$ small enough, we can achieve that
$\omega_{\eps}$ is so close to the vertical sides of $\partial R$,
that it is contained in the set where $a(x,y)$ is bounded away
from $0$. It follows that
\begeq
\label{eq18}
\norm{\chi u}_{L^2}^2 \leq {\cal O}(1)\norm{f}_{L^2}^2+{\cal
O}(1)\abs{\lambda}^2 \int_{\Omega} a(x,y)\abs{u(x,y)}^2\,dx\,dy.
\endeq

Now observe that multiplying (\ref{eq16}) by $\overline{u}$, we
obtain, integrating by parts and taking the imaginary part, that
\begeq
\label{eq19}
2 \lambda \int_{\Omega} a\abs{u}^2\,dx\,dy=\Im \int_{\Omega} f\overline{u}\,dx\,dy.
\endeq
Hence
\begeq
\label{eq20}
\abs{\lambda} \int_{\Omega} a\abs{u}^2\,dx\,dy\leq \norm{f}_{L^2}\norm{u}_{L^2}.
\endeq
and using this in (\ref{eq18}) we see that
\begeq
\label{eq21}
\norm{\chi u}_{L^2}^2 \leq {\cal O}(1)\norm{f}_{L^2}^2+{\cal
O}(1)\abs{\lambda}\norm{f}_{L^2}\norm{u}_{L^2}.
\endeq

It remains to estimate the $L^2$-norm of $(1-\chi)u$, and to that end we notice that the
support of this function is contained in the set where $a$ is bounded away from zero. Therefore, another application
of (\ref{eq20}) shows that
\begeq
\label{eq22}
\norm{(1-\chi)u}_{L^2}^2\leq {\cal O}(1) \norm{f}_{L^2}\norm{u}_{L^2}, \quad \abs{\lambda}\geq 1.
\endeq
Putting together the estimates (\ref{eq21}) and (\ref{eq22}) we
get
\begeq
\label{eq23}
\norm{u}_{L^2}\leq {\cal
O}(1)\norm{f}_{L^2}+{\cal
O}(1)\abs{\lambda}^{1/2}\norm{f}_{L^2}^{1/2}\norm{u}_{L^2}^{1/2},
\endeq
so that
\begeq
\label{eq24}
\norm{u}_{L^2}\leq {\cal O}(1)\abs{\lambda}\norm{f}_{L^2}.
\endeq
This completes the proof.
\end{proof}

\medskip
With Proposition 2.2 available, we are in the position to
estimate the norm of the resolvent of the operator ${\cal A}$ in (\ref{eq1.2}) on the real
axis. In doing so, we notice that a simple computation, as in~\cite{Lebeau} and~\cite{Hi}, shows that
\begeq
\label{eq24.5}
\left(\lambda-{\cal A}\right)^{-1}=\left(
\begin{array}{cc}
 R(\lambda)(2ia-\lambda) & -R(\lambda) \\
 R(\lambda)(2ia\lambda-\lambda^2)-1& -\lambda R(\lambda)
\end{array}\right),\quad \lambda \in \real.
\endeq
Here, as in Proposition 2.2, we have written $R(\lambda)=(-\Delta+2i\lambda a-\lambda^2)^{-1}$. To estimate the norm
$(\lambda-{\cal A})^{-1}$ as an operator on ${\cal H}=H^1_0\times L^2$,
we have to derive bounds on the operators
\begeq
\label{eq25}
R(\lambda): L^2\rightarrow H^1_0,\quad R(\lambda)(2ia\lambda-\lambda^2)-1: H^1_0\rightarrow L^2,
\endeq
and
\begeq
\label{eq26}
R(\lambda)(2ia-\lambda): H^1_0\rightarrow H^1_0.
\endeq

Now Proposition 2.2 together with an integration by parts argument shows that
\begeq
\label{eq26.5}
R(\lambda)={\cal O}(\lambda^2): L^2\rightarrow H^1_0,
\endeq
and hence by duality,
\begeq
\label{eq27}
R(\lambda)={\cal O}(\lambda^2): H^{-1}\rightarrow L^2.
\endeq
When estimating the norm of the second operator in (\ref{eq25}),
we use that $R(\lambda)(2ia\lambda-\lambda^2)-1=R(\lambda)\Delta$, and
hence, combining (\ref{eq27}) together with the fact that the Laplacian
$\Delta: H^1_0\rightarrow H^{-1}$ is continuous, we get
\begeq
\label{eq28}
R(\lambda)(2ia\lambda-\lambda^2)-1=R(\lambda)\Delta={\cal O}(\lambda^2): H^1_0 \rightarrow L^2.
\endeq

It remains to estimate the norm in (\ref{eq26}), and to that end, we write
\begeq
\label{eq28.5}
R(\lambda)(2ia-\lambda)=\frac{1}{\lambda}(1+R(\lambda)\Delta).
\endeq
If $f\in H^1_0$ and $u=R(\lambda)\Delta f\in H^1_0$ then
\begeq
\label{eq29}
\left(-\Delta+2ia\lambda-\lambda^2\right)u=\Delta f\in H^{-1},
\endeq
and hence multiplying (\ref{eq29}) by $\overline{u}$, integrating by parts, and taking the real part, we get
\begeq
\label{eq30}
\norm{u}_{H^1_0}^2-\lambda^2\norm{u}_{L^2}^2\leq
\norm{\Delta f}_{H^{-1}}\norm{u}_{H^1_0}\leq \norm{f}_{H^1_0}\norm{u}_{H^1_0}.
\endeq
Therefore,
\begeq
\label{eq31}
\norm{u}_{H^1_0}^2\leq {\cal O}(1)\left(\lambda^2 \norm{u}_{L^2}^2+\norm{f}_{H^1_0}^2\right).
\endeq
When estimating the $L^2$-norm of $u=R(\lambda)\Delta f$, we use (\ref{eq28}) to conclude that
\begeq
\label{eq32}
\norm{u}_{L^2}\leq {\cal O}(\lambda^2)\norm{f}_{H^1_0},
\endeq
and combining this estimate with (\ref{eq31}) we get
\begeq
\label{eq33}
\norm{u}_{H^1_0}\leq {\cal O}(\abs{\lambda}^3)\norm{f}_{H^1_0}.
\endeq
An application of (\ref{eq28.5}) then shows that
\begeq
\label{eq34}
R(\lambda)(2ia-\lambda)={\cal O}(\lambda^2): H^1_0 \rightarrow H^1_0.
\endeq

Combining Proposition 2.2 together with (\ref{eq26.5}), (\ref{eq28}),
(\ref{eq34}), and the fact that ${\cal A}$ has no real eigenvalues, we get the basic bound
\begeq
\label{eq35}
(\lambda-{\cal A})^{-1}={\cal O}((1+\abs{\lambda})^2): {\cal H}\rightarrow {\cal H},\,\, \lambda \in \real.
\endeq

\medskip
We shall finally show how the bound (\ref{eq35}) allows us to conclude
the proof of Theorem 1.1. In doing so, we shall follow the argument
of~\cite{Lebeau} closely, adapting it to the present case.

When $k>1$ is an integer and $x\in {\cal H}$, we write, as in~\cite{Lebeau}
and~\cite{Burq}, for $t>0$,
\begeq
\label{eq38}
e^{it{\cal A}}(1-i{\cal A})^{-k}x=\frac{1}{2\pi i}\int_{\gamma}
e^{it\lambda} \frac{1}{(1-i\lambda)^k}(\lambda-{\cal A})^{-1}x\,d\lambda.
\endeq
Here $\gamma=\{\lambda\in \comp; \lambda=\eta-i/2,\,\,\eta\in \real\}$. Furthermore,
when $X=X(t)=\gamma_1 (t/\log t)^{1/2}$, where $\gamma_1>0$ is to
be chosen, we use the same decomposition of (\ref{eq38}) as
in~\cite{Lebeau},
\begin{eqnarray}
\label{eq39}
& & e^{it {\cal A}}(1-i{\cal A})^{-k}x \\ \nonumber
& = & \frac{1}{2\pi
i}\frac{1}{\sqrt{2\pi}}\int_{\gamma}\int_{{\bf R}}
e^{it\lambda}\frac{1}{(1-i\lambda)^k}e^{-(\lambda-\tau)^2/2}(\lambda-{\cal
A})^{-1}x\,d\tau\,d\lambda \\ \nonumber
& = & \int_{\gamma}\int_{\abs{\tau }\leq
X}\cdots\,+\int_{\gamma}\int_{\abs{\tau }\geq X}\cdots=:I_1+I_2,
\end{eqnarray}
where
\begeq
\label{eq42}
I_1=\frac{1}{2\pi i}\frac{1}{\sqrt{2\pi}}\int_{\gamma}\int_{\abs{\tau }\leq X}
e^{it\lambda}\frac{1}{(1-i\lambda)^k}e^{-(\tau-\lambda)^2/2}(\lambda-{\cal
A})^{-1}x\,d\tau\,d\lambda,
\endeq
and
\begeq
\label{eq42.5}
I_2=\frac{1}{2\pi i}\frac{1}{\sqrt{2\pi}}\int_{\gamma}\int_{\abs{\tau }\geq X}
e^{it\lambda}\frac{1}{(1-i\lambda)^k}e^{-(\tau-\lambda)^2/2}(\lambda-{\cal
A})^{-1}x\,d\tau\,d\lambda.
\endeq

Now combining (\ref{eq35}) together with a perturbation argument, we see that
the function $\lambda \mapsto (\lambda-{\cal A})^{-1}x$ is holomorphic in $\lambda$
with values in ${\cal H}$, in the region below and including the curve
\begeq
\label{eq43}
\gamma_{\eps_0}=\left\{\lambda=\eta+i\frac{\eps_0}{\eta^2},\,\,\abs{\eta}\geq
1\right\}\cup \left\{ \lambda=\eta+i\eps_0,\,\,\abs{\eta}\leq 1\right\},
\endeq
where $\eps_0>0$ is small enough but fixed. Also, we see that along $\gamma_{\eps_0}$, $(\lambda-{\cal A})^{-1}x$
is bounded in ${\cal H}$ by ${\cal O}(1)(1+\abs{\eta})^2\norm{x}_{{\cal H}}$. Therefore, when estimating
$I_1$, we may write
\begeq
\label{eq44} I_1=\frac{1}{2\pi
i}\frac{1}{\sqrt{2\pi}}\int_{\gamma_{\eps_0}}\int_{\abs{\tau }\leq X}
e^{it\lambda}\frac{1}{(1-i\lambda)^k}e^{-(\tau-\lambda)^2/2}(\lambda-{\cal
A})^{-1}x\,d\tau\,d\lambda.
\endeq
The contribution to $I_1$ coming from the part of
$\gamma_{\eps_0}$ where $\abs{\eta}\leq 1$ is easily seen to be
bounded by ${\cal O}(e^{-\eps_0 t}X(t))\norm{x}_{{\cal H}}$, hence
decaying exponentially, and therefore we may concentrate on the
part of $\gamma_{\eps_0}$ where $\abs{\eta}\geq 1$. When
estimating the contribution coming from this part, exactly as
in~\cite{Lebeau}, we treat separately the cases when $1\leq
\abs{\eta}\leq \gamma_2 (t/\log{t})^{1/2}$ and $\abs{\eta}\geq
\gamma_2 (t/\log {t})^{1/2}$. Here $\gamma_2>0$ is to be chosen.
The ${\cal H}$--norm of the integrand in $I_1$ coming the part
where $\abs{\eta}\geq 1$ is bounded by
\begeq
\label{eq45}
{\cal O}(1)
\exp\left(-\frac{t\eps_0}{\abs{\eta}^2}\right)\frac{\abs{\eta}^2}{(1+\abs{\eta})^k}e^{-(\tau-\eta)^2/2}\norm{x}_{{\cal
H}}.
\endeq
Here $\abs{\tau}\leq \gamma_1 (t/\log{t})^{1/2}$. It follows that
the contribution to $I_1$ coming from the region where $1\leq
\abs{\eta}\leq \gamma_2 (t/\log {t})^{1/2}$ is controlled by
\begeq
\label{eq46}
{\cal O}(1)\frac{t}{\log
t}\frac{1}{t^{\eps_0/\gamma_2^2}}\norm{x}_{{\cal H}},
\endeq
which decays as any fixed inverse power of $t$, provided that
$\gamma_2>0$ is sufficiently small. Having fixed $\gamma_2>0$
small enough, as in~\cite{Lebeau}, we choose $\gamma_1\in
(0,\gamma_2)$, and observe that for $\abs{\tau}\leq \gamma_1
(t/\log {t})^{1/2}$ and $\abs{\eta}\geq
\gamma_2(t/\log{t})^{1/2}$, we have
\begeq
\label{eq47}
\frac{1}{2}(\tau-\eta)^2\geq \delta\left(\tau^2+\eta^2\right),\quad \delta>0.
\endeq
It follows then that the final contribution to $I_1$ coming from the region
where $\abs{\eta}\geq \gamma_2 (t/\log{t})^{1/2}$ is obtained by
integrating over this region
\begeq
\label{eq48}
{\cal O}(1)\left(\frac{t}{\log
t}\right)^{1/2}\frac{\abs{\eta}^2}{(1+\abs{\eta})^k} e^{-\delta \eta^2}\norm{x}_{{\cal H}},
\endeq
and this clearly decays rapidly (and even exponentially) as
$t\rightarrow \infty$. We conclude that for each $N\in \nat$,
\begeq
\label{eq49}
I_1={\cal O}_N(t^{-N})\norm{x}_{{\cal H}}.
\endeq

When estimating $I_2$ in (\ref{eq42.5}), we continue to follow~\cite{Lebeau} and write
\begeq
\label{eq50}
I_2=e^{it{\cal A}}J,
\endeq
where
\begeq
\label{eq51}
J=\frac{1}{2\pi i}\frac{1}{\sqrt{2\pi}}\int_{\abs{\tau }\geq X}\int_{\gamma}
\frac{1}{(1-i\lambda)^k}e^{-(\tau-\lambda)^2/2}(\lambda-{\cal
A})^{-1}x\,d\tau\,d\lambda.
\endeq
Since $e^{it{\cal A}}$ is uniformly bounded on ${\cal H}$ for
$t\geq 0$, it suffices to estimate $J$ in the energy norm, as
$t\rightarrow \infty$. Exactly as in~\cite{Lebeau}, we then see that we have to estimate
the integral
\begeq
\label{eq52}
\int_0^{\infty} S_{X-\frac{i}{2}+\mu e^{-i\pi/8}}(x)\,d\mu,
\endeq
where
\begeq
\label{eq53}
S_{\tau}(x)=\int_{\gamma} \frac{1}{(1-i\lambda)^k} e^{-(\lambda-\tau)^2}
\left(\lambda-{\cal A}\right)^{-1}x\,d\lambda.
\endeq
Arguing precisely as in~\cite{Lebeau} and making a contour deformation as in that paper, we then verify that with
$\tau=X-\frac{i}{2}+\mu e^{-i\pi/8}$, $\mu\geq 0$,
\begeq
\label{eq54}
\norm{S_{\tau}(x)}_{{\cal H}}\leq {\cal O}(1)\left(e^{-\delta(\mu^2+X^2)}+\frac{1}{(1+\mu)(\mu+X)^k}\right)
\norm{x}_{{\cal H}},\quad \delta>0.
\endeq
Let us remark here that when deriving (\ref{eq54}), the bound (\ref{eq1.5}) is important.

It follows from (\ref{eq50}), (\ref{eq52}), and (\ref{eq54}) that the ${\cal H}$-norm of $I_2$ does
not exceed a constant times the norm of the vector $x$ multiplied by
\begeq
\label{eq55}
\int_0^{\infty} \left(e^{(-\delta(\mu^2+X^2))}+\frac{1}{(1+\mu)(\mu+X)^k}\right)\,d\mu.
\endeq
Since
\begeq
\label{eq56}
\int_0^{\infty} \frac{d\,\mu}{(\mu+1)(\mu+X)^k}={\cal O}_k\left(\frac{\log X}{X^k}\right),
\endeq
we conclude, recalling the definition of $X$, that
\begeq
\label{eq57}
I_2={\cal O}_k(1)\frac{(\log t)^{k/2+1}}{{t}^{k/2}}\norm{x}_{{\cal H}},\quad t\geq 2,\quad k=2,3,\ldots
\endeq

This completes the proof of Theorem 1.1 in the case when $k>1$ is an integer. Using an interpolation argument
as explained in~\cite{Burq}, we get the result for a general $k>0$.

{\it Remark}. By refining again the analysis above as in~\cite[Section 4]{Burq}
one could probably avoid the logarithmic loss.

\section{Improved decay estimates: an example}\label{section3}
\setcounter{equation}{0}
In this section, we shall give an example of a class of damping
coefficients vanishing away from a \neigh{} of the non-rectangular
part of the domain, for which the result of Theorem 1.1 can be improved,
leading to an estimate of the type (\ref{eq3.6}). When doing so, to fix the ideas, we shall let $\Omega=S$ be the
Bunimovich stadium, defined in the example preceding the formulation of Theorem 1.1.

When $S=R\cup W$, $R=[0,1]_x\times [0,\pi]_y$, we let
$0\leq a\in C^{\infty}(\overline{S})$ be such that $a>0$ in
$\overline{W}$. We assume that
\begeq
\label{eq57.8}
a^{-1}(0)=[\delta, (1-\delta)]_x\times [0,\pi]_y,\quad 0<\delta \ll 1,
\endeq
and that $a$ is independent of $y$ when $x$ is close to $\delta$ and $1-\delta$.
Furthermore, let us assume for simplicity that close to $\delta$ (respectively $1-\delta$), we have
\begeq
\label{eq58}
\frac {d^m} {dx^m}a(x)\leq 0 \qquad (\text{respectively} \geq 0),
\endeq
for some $m\geq 4$. We then immediately deduce the following
result.
\begin{lemma}
For any $n<m$ there exists $C_{n,m}>0$ such that
$$
\abs{a^{(n)}(x)} \leq C_{n,m} a(x)^{\frac {m-n} {m}}.
$$
\end{lemma}
\begin{proof}
It suffices to consider the case $x\geq (1-\delta)$. Using Taylor's formula, we get for $x\geq (1- \delta)$,
$$
a(x) = \int_{1- \delta}^x \frac {(x-s)^{m-1}}{(m-1)!} a^{(m)}(s) ds, \qquad
a'(x) = \int_{1- \delta}^x \frac {(x-s)^{m-2}}{(m-2)!} a^{(m)}(s) ds.
$$
As a consequence,
\begeq
\label{eq58.1}
\abs{a'(x)} \leq \int_{(1- \delta)} ^y \frac {(x-s)^{m-2}}{(m-2)!} a^{(m)}(s)\,ds +
\int_{y} ^x\frac {(x-s)^{m-2}}{(m-2) !} a^{(m)}(s)\,ds,
\endeq
with $(1- \delta) < y<x$ to be chosen. To estimate the first integral in (\ref{eq58.1}) we use the bound
$$
\frac {(x-s)^{m-2}}{(m-2)!}= \frac{(m-1)}{{x-s}}\frac {(x-s)^{m-1}}{(m-1)!} \leq  \frac{(m-1)}{x- y}
\frac {(x-s)^{m-1}}{(m-1!},
$$
and for the second integral, we use that
$$
\frac {(x-s)^{m-2}}{(m-2)!}\leq \frac {(x-y)^{m-2}}{(m-2)!}.
$$
We obtain
$$
|a'(x)| \leq \frac{(m-1)}{(x-y)} a(x) + {\cal O}_m(1) (x-y)^{m-1}.
$$
Choosing $y$ so that  $a(x) = (x-y)^m$ gives the lemma for $n=1$. The general case is similar.
\end{proof}

Let us consider now the stationary problem
\begeq
\label{eq57.5} \left(-\Delta+2i\lambda
a-\lambda^2\right)u=f\in L^2(S),\quad u|_{\partial S}=0,\quad \lambda \gg 1.
\endeq
As in the proof of Proposition 2.2, we let $0\leq \chi_\lambda\in C^{\infty}_0((0,1))$ be a cut-off function
so that $\chi_\lambda$ vanishes for $x$ close to $0$ or $1$.Then the function $\chi_\lambda u$ vanishes on $\partial
R$ and satisfies in the interior of the rectangle,
\begeq
\label{eq58.5} \left(-\Delta-\lambda^2\right)\chi_\lambda u=\chi f+\chi_\lambda''
u-2\partial_x(\chi_\lambda' u)-2i\lambda a(x)\chi_\lambda u.
\endeq
We choose $\chi_\lambda=\chi (\lambda a(x))$ with $\chi=0$ for $\abs{x} \geq 2$, and $\chi=1$ for $\abs{x}\leq 1$,
so that in the support of $\chi_\lambda'(x)$ we have
\begeq
\label{eq59}
a(x)\sim \frac{1}{\lambda},\quad \lambda \gg 1.
\endeq
From Lemma 3.1 we get the following bounds on the derivative of $\chi_\lambda$,
\begeq
\label{eq61.1}
\abs{\chi_\lambda '} = \abs{\lambda \chi'(\lambda a(x)) a'(x)} \leq {\cal O}(1) \lambda^{ \frac 1 m},
\endeq
and similarly
\begeq
\label{eq61.1bis}
\abs{\chi_\lambda ^{(n)}} \leq {\cal O}(1) \lambda^{ \frac n m}, \quad n< m.
\endeq

We now come to estimate the $L^2$-norms of the functions $\chi_\lambda''
u$, $\chi_\lambda' u$, and $\lambda a(x)\chi_\lambda u$, occurring in the right
hand side of (\ref{eq58.5}). In doing so, we write
$$
\chi_\lambda' u=\frac{\chi' a^{1/2} u}{a^{1/2}},
$$
and using (\ref{eq59}), (\ref{eq61.1}), we get
$$
\abs{\frac{\chi_\lambda '} {a^{1/2}}} \leq {\cal O}(1) \lambda ^{\frac 1 2 + \frac 1 m},
$$
and consequently
\begeq \label{eq62} \norm{\chi_\lambda ' u}_{L^2(R)}\leq
{\cal O}(1) \lambda^{\frac 1 2 + \frac 1 m}
\norm{a^{1/2}u}_{L^2(S)}.
\endeq
Estimating the $L^2$-norm of $\chi_\lambda'' u$ in a similar way, we get
\begeq \label{eq62.5} \norm{\chi_\lambda'' u}_{L^2}\leq {\cal
O}(1)\lambda^{\frac 1 2 + \frac 2 m} \norm{a^{1/2}
u}_{L^2(S)}.
\endeq
Finally, a similar argument shows that \begeq \label{eq63}
\norm{a\lambda \chi_\lambda  u}_{L^2(R)}\leq {\cal O}(1)
\lambda^{1/2}\norm{a^{1/2}u}_{L^2(S)}.
\endeq

At this point we can apply Proposition 6.1 of~\cite{BZ1}
to (\ref{eq58.5}) directly, choosing $\omega_x$ there to be a
\neigh{} of the edges $x=0$ and $x=1$ where $a$ is bounded from below and consequently  $\chi_\lambda u$ vanishes.
Using (\ref{eq62}), (\ref{eq62.5}), and (\ref{eq63}), we get, with $\omega=\omega_x\times [0,\pi]_y$,
\begeq \label{eq63.5}
\norm{\chi_\lambda u}_{L^2(R)}^2\leq {\cal
O}(1)\left(\norm{f}_{L^2(R)}^2+\lambda^{1+ \frac 4 {m}}
\norm{a^{1/2}u}^2_{L^2(S)}\right).
\endeq
An application of (\ref{eq20}) in (\ref{eq63.5}) gives next
\begeq \label{eq64}
\norm{\chi_\lambda u}_{L^2}^2\leq {\cal O}(1)\left(\norm{f}_{L^2}^2+
\lambda^{\frac 4 {m}}\norm{f}_{L^2}\norm{u}_{L^2}\right).
\endeq
It remains to control the $L^2$-norm of $(1-\chi_\lambda)u$, and when
doing so we remark that the support of $1-\chi_\lambda$ is contained in
the set where $a\geq 1/\lambda$, $\lambda \gg 1$. Using (\ref{eq20})
for the second time, we infer that
\begeq
\label{eq67} \norm{(1-\chi_\lambda)u}_{L^2}^2\leq {\cal
O}(1)\norm{f}_{L^2}\norm{u}_{L^2}.
\endeq
Putting together (\ref{eq64}) and (\ref{eq67}) we get the estimate
\begeq \label{eq68} \norm{u}_{L^2}\leq {\cal
O}(1)\left(\norm{f}_{L^2}+\lambda^{\frac 2 m}
\norm{f}_{L^2}^{1/2}\norm{u}_{L^2}^{1/2}\right),
\endeq
and finally,
\begeq \label{eq69} \norm{u}_{L^2}\leq {\cal
O}(1)\abs{\lambda}^{\frac 4 {m}} \norm{f}_{L^2},\quad \lambda\in
\real,\,\, \abs{\lambda}\gg 1.
\endeq

Repeating the arguments of section 2, with the bound (\ref{eq69}) in place of Proposition 2.2, we get the bound
\begeq
\label{eq70}
\left(\lambda-{\cal A}\right)^{-1}={\cal O}(1)\left(1+\abs{\lambda}\right)^{1+\frac 4{m}}:
{\cal H}\rightarrow {\cal H}, \quad \lambda \in \real,
\endeq
and Theorem 1.2 follows.

\end{document}